# The Moravian Crossroad: Mathematics and Mathematicians in Brno Between German Traditions and Czech Hopes

Laurent Mazliak[1] and Pavel Šišma[2]

**Abstract**. In this paper, we study the situation of the mathematical community in Brno, the main city of Moravia, between 1900 and 1930. During this time, the First World War and, as one of its consequences, the creation of the independent state of Czechoslovakia, led to the reorganization of the community. German and Czech mathematicians struggled to maintain forms of cohabitation, as political power slid from the Germans to the Czechs. We show how the most active site of mathematical activity in Brno shifted from the German Technical University before the war to the newly established Masaryk University after.

Following the Armistice in November 1918, national borders were drastically redrawn in Europe. New states, such as Poland, Czechoslovakia, and Yugoslavia, emerged and immediately had to situate themselves geopolitically. To the new states, military circumstances dictated a system of alliances with the victors of World War I—mainly Great Britain, France, Italy, and the United States. The Allies having imposed an embargo on scientific collaboration with Germany and Austria, mathematicians in the new states were even more inclined to translate the geopolitical situation into intellectual cooperation.[3] For countries that had hitherto been submitted to German, Austrian, or Russian influence and authority, this led to major cultural realignments. At the local levels, the tense cohabitation between different linguistic, religious, and national communities complicated matters even more. But what effect did this political reshuffling of territories have on mathematics, if any? One way to start examining this question is to look at correspondences that were established at that time between mathematicians in Allied nations and those who lived in the newly created states of Central Europe.

In 1919, the 40-year old mathematician Maurice Fréchet (1878–1973) was put in charge of organizing the mathematics department at the University of Strasbourg that was reestablished as the intellectual showcase of French science in the buildings of the German Kaiser-Wilhelm University. On June 29, he addressed a letter to a colleague of his in Prague, the capital of Czechoslovakia, which had just declared its independence on October 29, 1918:

> Would you allow me to inquire about the universities that will be staying or created on the territory of your new state. Perhaps one of your students moreover would oblige me by sending the list of mathematics professors in Czechoslovakian universities, as well as the list of Czechoslovakian journals publishing original mathematical

---

[1]    Laboratoire de Probabilités et Modèles Aléatoires, Université Pierre et Marie Curie, Paris, France. laurent.mazliak@upmc.fr
[2]    Katedra matematiky Přírodovědecké fakulty Masarykovy univerzity, Brno, Czech Republic. sisma@math.muni.cz

[3] On this question see (among others) [Lehto 1998], [Mc Millan 2001] and [Mazliak & Tazzioli 2009].

> articles written by your fellow countrymen. Are any of these journals publishing articles in French?[4]

To whom Fréchet had written is not clear. What we know however is that, some months later, Bohuslav Hostinský (1884–1951), a mathematician and physicist then serving as secretary of Czechoslovakia's National Provisory Committee, was given Fréchet's letter. His good knowledge of French was the probable reason for his ending up with the letter [Bru 2003]. The son of the musicologist Otakar Hostinský, a prominent member of the Czech intelligentsia, Bohuslav Hostinský had defended a Ph.D.thesis on Lie spherical geometry in 1906 and spent the academic year 1908–1909 in Paris where he heard the lectures of Emile Picard, Henri Poincaré, and Gaston Darboux.

On October 19, 1919, Hostinský replied to Fréchet, informing the Strasbourg mathematician that new universities were to be opened in Bratislava and in Brno—the latter to be named Masaryk University after the first president of Czechoslovakia and where he was soon to be appointed as professor of theoretical physics. Hostinský mentioned that both of the major Czech journals, the Casopis pro pěstování matematiky a fysiky (Journal for the Cultivation of Mathematics and Physics) and the Véstník Královské české společnosti nauk (Bulletin of the Royal Czech Science Society), were about to change their language policy to increase the French and English presence at the expense of the German. To conclude his letter, Hostinský suggested to Fréchet that he be his main contact in Czechoslovakia in case he needed one.

This was the seed of an extensive, mathematically fruitful correspondence between the two mathematicians [Havlová et al. 2005]. Started in the aftermath of WWI, this correspondence, we argue, is representative of the way in which the emergence of independent states allowed the development of new cultural and scientific networks. In mathematics, this is particularly significant, as it closed a period of French and German domination. The example of Czechoslovakia is enlightening in that it shows how a national community developed its own local network of scientific institutions and how, while it has hitherto belonged to the German cultural sphere, this community started to look for more extensive contacts with other spheres. The same phenomenon can be found in different countries of Eastern and Central Europe : [Mazliak, 2012] provides an example drawn from the Bulgarian situation. The impressive development of mathematical research between the two world wars was particularly obvious in these countries where entirely new institutions or journals were created that allowed them to rival the best foreign mathematical research in a few years time[5].

In what follows, we focus on mathematical life in Brno, the capital of Moravia, before, during, and after WWI.[6] Bordering Austria on the South, Moravia has always been a multicultural region. Due to the presence of a strong German minority, it was one of the parts of Europe where the question of nationalities would be posed with special acuteness in 1918 (see figure 1). This complex cultural cohabitation was reflected in the fact that two distinct mathematical communities lived side by side in the town. The discrepant effects of WWI on both communities are examined below from an international standpoint.[7]

---

[4] "Voulez-vous me permettre de vous demander quelles sont les différentes universités qui doivent subsister ou être créées sur le territoire de votre nouvel Etat. Peut-être en outre un de vos étudiants voudrait bien me rendre le service de me communiquer la liste des professeurs de mathématiques des universités tchécoslovaques, et celles des périodiques tchécoslovaques qui impriment les mémoires originaux de math. par vos savants compatriotes. Y a-t-il un de ces périodiques qui soit publié en français?" [Hostinský's Papers].
[5] See for instance what Denjoy writes in [Denjoy 1939, p.2].
[6] To simplify, we shall always refer to the city by its Czech name, although when the German community is concerned it obviously should be replaced by Brünn. Similarly, we will throughout

## 1. The Czech Fight for Higher Education in Moravia Before WWI

Moravia is an interesting place to study the effect of the postwar cultural reconfigurations on mathematics. The father of Czechoslovakian independence, Tomáš G. Masaryk wrote that in the "so-called 'German territories' in Bohemia (Moravia and Silesia) numerous Czechs are living; it is therefore fair that the Czech state should keep them; it would be unfair to sacrifice several hundreds of thousand Czechs to the furor teutonicus" [Masaryk 1920]. This annexion was of course the seed of many conflicts to come. As Edvard Beneš, Czechoslovakia's Minister of Foreign Affairs, explained to the delegates of the Peace Conference in Paris, on February 5, 1919: "the relations of Czechoslovakia with its neighbours have to be settled in order to avoid any future conflict."[8]

In Brno the coexistence of the Czech majority with an important German minority was indeed tense. This fact is crucial in order to understand the shape of its educational institutions between 1880 and 1930. Although the German minority that had exerted a dominant role in the cultural arena suddenly lost its preeminent place in 1918, its cultural influence remained important. It is in this context that one must interpret efforts at building cultural bridges with Allied powers. Our study emphasizes a fact that may seem obvious at first glance, namely, that the history of mathematics in Moravia cannot avoid taking into account the relationship between communities, even if the contours of these communities were never very precise. These contacts were complicated, mixing rivalry and dialogue, and most often reduced to the minimum.[9]

---

write Strasbourg where the Germans would of course have written Straßburg. We will not, however, use today's denomination for university towns of the Austro-Hungarian Empire, since the context of our study is clearly German: accordingly, Chernivtsi will be called Czernowitz, Lviv Lemberg, and Wrocław Breslau.

[7] Who, in this story, were the "Germans" and who were the "Czechs" is an important but difficult point to consider. During the period we study, the answer to such question was never univocal and was constantly changing depending on political and social conditions. In his study of the German minority in Czechoslovakia, Křen has shown the fluctuation in the definition of Germans (and Czechs) [Křen 1998]. In population censuses taken in the years 1880–1900, the numbers of inhabitants of Brno who declared Czech as their main language of communication, varied from 30 to 40 %. In the last census before WWI (1910), 41,000 out of 126,000 Brno inhabitants (32%) declared Czech as their usual language [Dřímal & Peša 1973, vol. 2, p. 64]. These figures must however be considered with care. Political and economical reasons probably led to an overestimation of the German settlement. This can be inferred from the fact that in the very different political context of January 1919, 61% of the (almost identical) population declared that it belonged to the Czech community.

[8] Quoted from the newspaper *Le Matin*, February 6, 1919.

[9] One may here recall the well-known fact that the coexistence of several cultural communities in Brno of course ended with violence. The German invasion of 1939 was followed by the terrible years of occupation and the general expulsion of German-speaking inhabitants between 1946 and 1948.

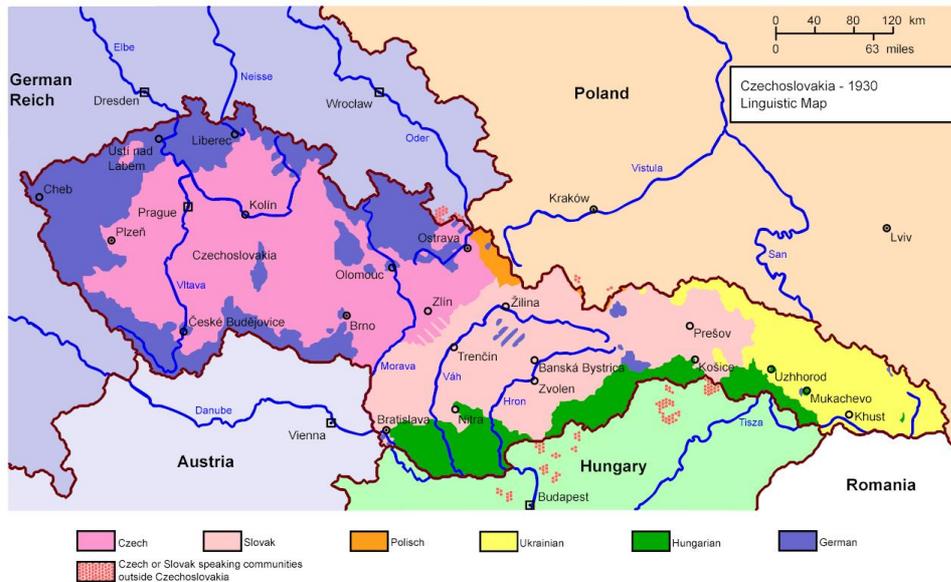

Figure 1. Czechoslovakian Linguistic Map in 1930. As can be seen, Brno is located on the border between the Czech and German zones. By derivative work: Henry Mühlpfordt, Mariusz Pazdziora (Czechoslovakia1930linguistic.jpg) [CC-BY-3.0 (http://creativecommons.org/licenses/by/3.0)], via Wikimedia Commons.

The first university to be established in Brno was, characteristically, a German institution. Founded in 1873, the Technical University (Technische Hochschule) replaced the Polytechnicum established in 1849, itself a distant heir of the old Olomouc Academy.[10] The Technical University was divided into faculties and was managed by an elected rector. Though the number of professors increased, the number of students stagnated and the Brno Technical University was, in fact, a small institution [Hellmer 1899]. But this does not mean that positions there were unappealing. Numerous Austrian scientists began their academic careers in modest size institutions. As Havránek recalls, the epigrammatic characterization of the professor's career in the Habsburg monarchy might have been: "Sentenced to Czernowitz, pardoned to Graz, promoted to Vienna" [Havránek 1998, p. 216]. Albeit less prestigious than Graz because of what was perceived as the hostile Czech environment, the Brno Technical University was certainly attractive due to its proximity to Vienna.

Since the beginning of the Czech national revival at the end of the 19th century, communities opposed one another on symbolic grounds and in cultural and intellectual life. In the 1860s, students from the Prague Technical University demanded to be and were lectured in Czech. In 1869, the Prague Polytechnicum was divided into two separate institutions and so was the venerable Charles University of Prague in 1882. Mathematicians were active in this movement. In 1862, a Czech counterpart to the German Mathematicians' Union [Deutsche Mathematiker-Vereinigung, hereinafter DMV] was founded by students in Prague: the Union of Czech Mathematicians and Physicists (Jednota českých matematiků a fyziků). Although the Jednota was initially devoted to the improvement of the students' scientific knowledge and lecturing skills without consideration of language or nationality, activists rapidly transformed it into a Czech national organization, isolating Czech from German students and severing ties with German teachers. The Jednota stirred the Czech intelligentsia's national consciousness. Although members of this society were spread over Bohemia and Moravia, official meetings and lectures were all held in Prague until WWI.[11]

---

[10] For an overview of the history of Czech mathematics, see [Nov´y 1996]; see also [Vetter 1923].

Until the 1860s, mathematicians at universities in Prague and Olomouc wrote their theses and scientific production in German whether Czech or German themselves. During the national revival, the situation changed: mathematicians split along national affiliations while career prospects started to diverge [Nový 1961, p. 221]. German mathematicians from Czech lands maintained tight contacts with mathematicians in Germany and Austria and were considered German mathematicians in Europe.[12] They hoped for and often got appointments to more prominent academic centers in the German-speaking world. More isolated, Czech mathematicians had, on the other hand, fewer opportunities for pursuing research. It is not so much that they lacked contact with European mathematicians before the war. After graduation, they often spent a year abroad, in Germany, France or Italy. But, at the end of the 19th century, Czech mathematicians had very small chances of being appointed in Austria or Germany. At most they could hope for a position in one of the two Czech universities in Prague. Back home, they were moreover asked to write textbooks in Czech, not only for university students, but also for secondary school pupils, while German mathematicians wrote textbooks on advanced topics.

Although present everywhere, this conflict took different forms in Bohemia and in Moravia. In Prague, due to the large Czech population, the Germans often experienced cohabitation as a threat. The division of the Prague university in 1882 into different German and Czech entities was a case in point. In 1897 German students, allying themselves with northern and western Bohemian Germans belonging to the movement Los von Prag (Out from Prague), tried to move German universities from Prague to Reichenberg (Liberec) [Cohen 2006]. Jaroslav Hašek's famous book *The Good Soldier Svejk* [Hašek 1930], in which Austrian militarism and bureaucracy are ridiculed, certainly exaggerated the tension in the 20th century. Still, heated debates sometimes pitted academics against each other. In 1907, mocking the requests made by the Prague Czech University to obtain more substantial subsidies, the rector of the (German) Charles University August Sauer for example called the university "a spoiled child," to which the politically-active physician Otakar Srdínko replied in a small brochure [Srdínko 1908].

Nearer to the border, South Moravia where the German population constituted a significant minority experienced cohabitation differently. Before the foundation of universities in Brno, Czech Moravian students were mostly attracted by Prague and the Germans by Vienna. At the end of the 1880s, out of roughly 1000 Moravians who received a higher education, about 700 went to Vienna, 250 to the Czech

---

[11] For details on the Union of Czech mathematicians and physicists, see [Seidlerová 1998].
[12] As Seidlerová observes, the situation continued after the independence of Czechoslovakia and even inside the country. She writes: "Even Czech university teachers and researchers often had no idea that in their works they actually cited a colleague from the Brno German Technical University" [Seidlerová 1998, p. 234]. We can also mention the case of Friedrich Urban, a German mathematician in Brno who, in 1923, wrote a book about random events in chains (in modern terms, Markov chains) precisely at the time when Hostinský started to be interested in them. Although in Brno at the time, Hostinský learned of it through Fréchet [Havlová et al. 2005].

University in Prague and 60 to the German University in Prague. By the mid-1890s, the number of Moravian students had increased to some 1300 and the question of opening a new university in Brno was raised. That in the Hasburg monarchy some towns with much smaller pools of students such as Czernowitz [today, Tchernivtsi] in distant Bukovina possessed a university where less than 300 students were enrolled, was seen as further evidence for the need of establishing a university in Brno [Jordán 1969, p. 40 & p. 43].

The central problem was the teaching language at the new university. Whereas a Czech population of 5 million had only one university to turn to in Prague, 8 million Germans could chose among five universities (Vienna, Prague, Graz, Innsbruck, and Czernowitz). In 1885, noticing that Czechs from Moravia and Silesia formed a sizeable proportion of students enrolled in Austrian universities, Masaryk pleaded for the establishment of a second Czech university [Masaryk 1885, p. 275]. His goal, as he explained later, was to insure a healthy competition and new positions for young Czech professors [Masaryk 1894]. If all young and talented graduates were to be appointed at the only Czech university that existed, then none other would stand a chance in future decades. By the late 1890s, Masaryk's fear was realized and a large number of privatdocents at the Czech University in Prague were faced with no professional prospects for promotion to a professorship.

In 1896, the Imperial Government conceded the right to establish a university as well as a new Czech technical university in Moravia. But this seemed like nothing but a formal declaration, as Vienna simultaneously asked both national components first to reach a preliminary agreement about the creation of a new university [Jordán 1969, pp. 47–48]. In 1899, a Czech technical university was opened in Brno as a counterpart to the German, but projects for establishing a Czech university were met with strong resistance. When Masaryk renewed the call for a Czech university in Brno, prominent members of the German community voiced their worries, especially since they felt that it was the (richer) German community that would have to foot the bill with its own taxes.[13] They warned that a Czech university would necessarily sink to a low academic standing, since language was a barrier to any attempt at establishing and maintaining contacts with world-leading scientists [Bachmann 1902]. To such fears of "Czechization," Sauer replied:

> Where will the students of the University of Brünn . . . find employment if not by occupying a position that until today has stayed in the hands of other races, and above all the Germans? . . . If Czechization means that Germans will be expelled from their positions, that Germans will be dominated, that Germans will be oppressed [by the Czechs], then this proud slogan must not be used to establish a second Czech university without our strongest protest against this villainy (quoted in [Srdínko 1908, p. 8]).

And so was the matter of a Czech university in Brno put to rest until the outbreak of World War I, despite a new series of petitions presented by Masaryk in 1912 to the Austrian Parliament [Jordán 1969, p. 105].

---

[13] To this claim, Srdínko replied with anger: "He who wants to turn civilization into an object of trade and offer it only to the rich, behaves as a barbarian" [Srdínko 1908, p. 9].

## 2. Mathematics in Brno Before 1914

Consider now the state of mathematics in higher education in Brno between 1880 and 1914. As we have argued, the split between communities makes it more appropriate to consider German and Czech institutions successively.

2.1. German Institutions. At the German Technical University, the chairs (Lehrstühle) of mathematics, of descriptive geometry, and of theoretical mechanics changed hands more frequently than in other departments. This may have been the case mainly because mathematicians had more opportunities to find other position, as mathematics was taught both in technical universities and in universities. Mathematicians mostly came to Brno from Graz, Prague, or Vienna where they were born, had relatives and had studied at local universities. Many of them, therefore, used the first opportunity to leave Brno if a chance to go back home presented itself [Šišma 2002].

The first renowned mathematical professor at the Brno Technical University was Emanuel Czuber (1851–1925) who came from the Prague German Technical University in 1886. A Czech from Prague, his original name, Cubr, had been Germanized during his studies at the German Technical University in Prague. He would soon leave Brno for the Vienna Technical University—a fine example of the way in which Brno seemed at the time destined to serve as a springboard to Vienna [Doležal 1928].

With Georg Hamel's (1877–1954) appointment in 1905, a brilliant period for mathematics began at the Brno German Technical University. The 28-year-old mathematician had defended a doctoral thesis in Göttingen about Hilbert's fourth problem and then a Habilitation in Karlsruhe. Hamel's presence and unprecedented activity, in particular in recruiting first-rate collaborators, seems to have given a boost to Brno's mathematical life. The Brno seminar was announced in the Jahresbericht der Deutschen Mathematiker-Vereinigung, the journal of the German Mathematicians' Union (DMV). In the fall of 1905, a Society for Mathematics and Physics (Mathematisch-Physikalische Gesellschaft) was created, an event announced in the JDMV in 1907. Over the next few years the mathematical seminar met every fortnight.[14] The full list of speakers between 1906 and 1913 includes besides Hamel himself and his colleague in Brno, Emil Waelsch (1863–1927), Max Ungar, Richard von Mises (in Brno, 1907–1909), Emil Fischer (in Brno, 1902–1910), Ernst Fanta (an actuary from Vienna who gave weekly lecture on his trade in Brno), Hans Hahn, Arthur Haas, Alfred Haar, Gottfried Rückle, Paul Ehrenfest, Heinrich Tietze (in Brno, 1910–1919), Johann Radon. . . . If one compares then with the similar mathematical seminar at the mathematical society of nearby Vienna, one finds that speakers often overlapped in a given year. Tietze, Lothar Schrutka, von Mises and Hamel went several times to Vienna while they held positions in Brno. But while the Vienna seminar concentrated on strictly mathematical topics, Brno was slightly slanted towards mechanics and mathematical physics (with talks about Planck's results by Hamel, hydrodynamics by von Mises, electromagnetism and gravitation by Gustav Jaumann). Personal relationships certainly played a great role in the Brno seminar. Several mathematicians on the list had met as students in Göttingen (Hamel, Fanta and Haar, for example). Hahn, Tietze, Ehrenfest and Gustav Herglotz had formed a small inseparable group of students at the Vienna University.[15] After Hamel's departure for Aachen in 1912, the life of the Gesellschaft, though less active, was maintained until the beginning of WWI.

---

[14] See Jahresbericht der deutschen Mathematiker-Vereinigung, vol. 16, 1907, pp. 396–397; vol. 18, 1909, pp. 104–105; vol. 21, 1912, pp. 58–59; vol. 23, 1914, pp. 52–53.

2.2. Czech Institutions. At the time of its foundation in 1899, several mathematicians were appointed to the Brno Czech Technical University. Among the first members of the faculty we find Karel Zahradník (1848–1916), who became the first rector, Jan Sobotka (1862–1931), the professor of descriptive geometry, and Antonín Sucharda (1854–1907). All had studied in Prague and worked as secondary school teachers at some point in their career. Zahradník came from the University of Zagreb in Croatia where he had remained for more than 20 years. Sobotka and Sucharda, who were younger, had nevertheless travelled, respectively, to Vienna, Zurich, and Breslau, and to Göttingen, Munich, Paris, and Strasbourg. The professor of mechanics was Václav Řehořovský (1849–1911) who had a similar background. As we can see, the foundation of the Czech Technical University in Brno created four new professorships for the Czech mathematical community.

In the following decade, Sucharda retired and Sobotka left after he obtained a professorship at the Prague Czech University where he ultimately trained the majority of Czech geometers in the first half of 20th century [Urban & Vančura 1962]. In 1906, Matyáš Lerch (1860–1922) was appointed as Sucharda's replacement. In 1908, Miloslav Pelíšek (1855–1940) was named professor of descriptive geometry in Sobotka's place [Franěk 1969, vol. 1, pp. 233–240]. A brilliant number theorist, Lerch had been a tenured professor at the University of Fribourg, Switzerland, since 1896. Lerch had studied mathematics in Prague where he worked as Eduard Weyr's assistant in 1885. He habilitated at the Prague Technical University in 1886 and over the next ten years published more than one hundred mathematical articles [Frank 1953]. While in Fribourg Lerch had tried, without success, to get a position at the Prague Czech University after the deaths of František Josef Studnička (1836–1903) and Weyr. In 1906, Lerch was finally appointed professor of mathematics at Brno. According to his assistant and successor, Karel Čupr, Lerch's lectures had not changed since Fribourg. Technical University engineering students found them more suitable for mathematics students and protested publicly.

The foundation of the second Czech Technical University brought changes to the organization of the Union of Czech Mathematicians and Physicists (Jednota). In 1913, a section of the Jednota was officially founded in Brno. In fact, lectures and meetings had begun to be organized in Brno immediately after the establishment of the Czech Technical University. During the years 1901–1911, about 55 lectures were held there. Auditors were mostly recruited among Brno secondary school mathematics teachers [Košťál 1967, pp. 18–21]. But contact with the German Mathematical Society in Brno was minimal and no mathematician or physicist from the Czech Technical University ever gave a lecture in Hamel's seminar. Indeed, mathematical life in Brno before WWI, such as it was, seemed completely split along the lines of national and linguistic communities. We were unable to find a single instance of mathematical cooperation between a German and a Czech university at an official level (such as an exchange of professors, common lectures and so on) or even between two individual members of these universities!

---

[15] Details about Herglotz and Ehrenfest can be found in an interesting article [Huijnen & Kox 2007] describing the context of the talk delivered by Ehrenfest in Brno in February 1912, when, after having decided to leave St. Petersburg, he was on the lookout for a permanent position in another European university.

**3. Mathematics in Brno during and after WWI**

3.1. The War. At the outbreak of WWI, mathematical life in Brno came to a standstill. The activity of the German Mathematical Society stopped immediately. After Hamel's departure from Brno, the professorship of theoretical mechanics remained vacant until 1916. From 1908 on, Hamel had tried to leave for Germany mostly for personal reasons due to his marriage in Cologne in 1909. He first tried to obtain a chair in Hannover but it was given to Constantin Caratheodory. Hamel eventually obtained a position in Germany in Aachen in 1912. Hamel's departure from Brno deprived mathematical life in the German Technical University from a major source of energy and the outbreak of WWI gave a deathblow to it. Two professors of mathematics—Lothar Schrutka (1881–1945) and Heinrich Tietze (1880–1964)—served in the Army, Tietze as officer on the front and Schrutka as a teacher at a secondary military school in Vienna. More generally, in the German Technical University, 5 professors among 33 and 7 privatdocents over 23 went into the army, as well as 40 assistants and 34 other employees.[16] All mathematical lectures and lessons were delivered by the professor of descriptive geometry Emil Waelsch. His assistants Wilhelm Schmid (1888–1963) and Rudolf Kreutzinger (1886–1959) were both drafted and replaced by students. Waelsch's assistants were taken prisoner in Russia and were unable to get back to to Brno until 1920 and 1921, respectively. This created problems for the teaching of descriptive geometry after the war, a time when the number of students increased. During the war, however, student enrollment at the German Technical University decreased from 950 to roughly 100–200.

Surprisingly, courses at the Czech Technical University hardly suffered in the first years of the war, since all professors of mathematics and descriptive geometry (Zahradník, Lerch, and Pelíšek) were old men at the time and therefore not drafted. After Zahradník's death in 1916, his successor, Jan Vojtěch (1879–1953), was appointed, but only in 1918. Of 177 teachers and employees of the Czech Technical University in Brno, 64 persons served in the army [Franěk 1969, vol. 1, p. 103]. The number of regular students at the school decreased from 571 in the academic year 1913–1914, to 254 in 1914–1915, and finally to 90 students in 1916–1917.[17] The remaining students of these empty years were often younger than before the war and belonged to the classes not yet called to the army. In 1913–1914, students under 19 years old represented 30% of the total enrollment, while in 1914–1915 this number had risen to 49%. In 1917, some soldiers were allowed to come back to universities, and student enrollment rose to 368, many of whom were students having passed their first-year examinations in 1914–15, now enrolled in their second year after having spent up to two years in the ranks.

---

[16] These data were obtained from our study of the programs of the German Technical University in Brno for the years 1913–1919 and from [Haussner 1924, pp. 33–34].

[17] See lists of students for academic years 1913–1914 to 1919–1920 in [CTU Archives].

At the German and Czech Technical Universities, buildings were turned into military hospitals, which complicated the resumption of teaching in 1918 [Šišma 2002, pp. 146–148]. Both universities saw their financial situation deteriorate and one can estimate that the endowment for the German Technical University as well as for the Czech Technical University had been reduced by 40% during the war ([Haussner 1924, p. 31], [Franěk 1969, p. 104]).

3.2. Czechoslovakia and the German Technical University. The Czechoslovakian declaration of independence on October 29, 1918 brought important changes to university life where, as we have seen, German institutions had played such an important role. The situation of the German communities within the borders of the Czech lands became unstable. While the German deputies from Bohemia and Moravia petitioned in Vienna for the annexation of German-populated regions by Germany (north Bohemia and Silesia) or by Austria (South Bohemia and South Moravia), the Czechoslovakian Government started to organize the military occupation of the lands as early as November 1918. In March 1919, a violent repression against the Sudete Germans temporarily settled the question and the German in-habitants were forcibly incorporated into the Czechoslovakian state [Bělina et al. 1993].

In 1919, Czech nationalists spent a lot of energy to prove that though the German Technical University had been favored by Vienna during the Austrian domination, the new Czechoslovak Government had not decided to destroy it in retaliation. In Parliament, František Mareš (1857–1942), a professor of physiology and philosophy, who became rector of the Charles University in 1920–21, mentioned that, while the budget of the German Technical University was 707,000 crowns at the time of the Austrian Emprire, it had risen to 1,753,000 crowns by 1919.[18] In fact, the extreme fluctuation of exchange rates and inflation during and after the war make it difficult to appreciate the significance of these numbers.[19] The estimated inflation index rate for the crown in the Czech lands in October 1918 was 1876 (taking 100 in 1914 as a basis) [Šedivý 2001, p. 245].

In Prague as well as in Brno, German university leaders lost their higher positions and often uttered alarming declarations. As Germans boycotted the Czechoslovak National Assembly of 1919, the Czech majority was free to exert tighter control on German higher education. In the first election held in 1920, there were 72 German deputies elected (that is, deputies whose party's name included the word "German") out of 300. When the organization of higher education in the new state was discussed in the Czechoslovakian Parliament, a German professor from Charles University in Prague expressed his helplessness: "We are in an appalling situation: a great part of our university will be thrown out onto the street.

---

[18] Czechoslovakian Parliament Debates, 34thSession, February 27, 1919.
[19] One reason for hyperinflation was that in 1914, in order to make the war popular, the Austro-Hungarian Government had decided to pay double the price for the main articles of necessity (such as grain, cattle, and horses): the enormous amount of 30 billions crowns had been printed by the Austro-Hungarian bank, resulting in huge inflation [Rašín 1923, p. 23].

The situation is distressing and is best expressed by the words: homeless, without means, without rights."[20] In the parliamentary debate, on February 27, 1919, the physician Srdínko—who had meanwhile been elected as deputy—contested the honesty of such declarations and petitioned for a redistribution of public resources among Czech and German universities. He recalled that, even before the war, texts had been published by foreign authors condemning the disproportion of means between universities. He quoted an article published in the Revue générale from 1911 stating that "a brutal and obvious fact appears from this amount of documents. It is the extraordinary disproportion existing between the credits attributed to the German University and those attributed to the Czech University, if one takes into account the respective populations" [De Visscher 1911, pp. 102–103].[21]

Perhaps understandably, professors at the German Technical University in Brno felt much looser connections with the city and the new state than their Czech colleagues. Many had been trained in other Austrian universities. Before WWI they had seen little difference between Brno, Graz, or Innsbruck; for one thing, Brno was certainly more attractive than distant Lemberg or Czernowitz. But in 1919, the German Technical University in Brno was an institution whose main constituency was a minority that had been deprived of much of its prior political and cultural power. Even if demobilization caused student enrollment to rise in 1920–1921 to twice its prewar value, the future of the Technical University was unclear. Indeed the existence of two German technical universities in Prague and in Brno for the German minority became a political issue when some drew attention to the fact that in the new state three million Germans enjoyed the same number of technical universities as nine million Czechs and Slovaks. The Czechoslovakian Government however chose to avoid direct confrontation with the Germans on this point and maintained the status quo. In March 1919, the professors of the Brno German Technical University took an oath of loyalty to the Czechoslovakian Republic [GTU Archives, B34, 416].

At the time, Austrian universities worried about the fate of their "sisterinstitutions" in the Czech lands. On December 14, 1918, at a meeting of Austrian universities in Vienna, the possibility of a common future was discussed. Fearing that it would be impossible to pursue their activity in the Czech capital, the representatives of the Prague German University suggested that their university be moved to a town in Bohemia where the Germans would form the majority.[22] The Association of Austrian German Engineers made up more specific plans and suggested to transfer the Prague German Technical University to Usti-nad-Labem in northern Bohemia and officially transform the German Technical University in Brno into a branch of the Vienna Technical University.[23] On December 23, 1918 the Academic Senate of the Prague German University declared that if the German border regions of

---

[20] Quoted from Czechoslovakian Parliament Debates,34thSession, February 27,1919.
[21] This paper offers an interesting insight into academic life in the Czech lands on the eve of First World War, mostly in Prague, though several comments are made on the Moravian situation. The paper contains a lot of figures for both the Czech and the German universities. Interestingly, de Visscher writes : "Contraints de s'adresser à l'étranger, les Tchèques ont craint la domination du génie allemand. Ils cherchent à des sources moins envahissantes, plus lointaines et plus calmes, le soutien d'une tradition scientifique et d'utiles sympathies" [De Visscher 1911, p.109].
[22] Czechoslovakian Parliament Debates, 34th Session, February 27,1919.
[23] Lidové Noviny, Brno,January 19,1919.

12Czechoslovakia were reunited with Germany or Austria, the university should also be transferred. This declaration infuriated Mareš:

> "What does that mean? It means that the Academic Senate of the German University in Prague does not acknowledge the legality of this country [and acts] as if the German University in Prague was not a property of the Czechoslovakian Republic."[24]

In this delicate political context, the local mathematical community at the Brno German Technical University underwent many changes. In 1919, mathematics professor Tietze left for Erlangen and Fanta stopped giving his weekly lectures on actuarial mathematics in Brno. It is possible that other docents from Vienna, who had worked at the Brno Technical University before the war, also stopped commuting.

The difficulties faced by the German Technical University in filling the position left open by Tietze's departure is illustrative of the complexity of the new situation. The position of mathematics professor indeed remained vacant until 1923, when Karl Mayr (1884–1940) was appointed. Eight mathematicians applied when the position was first opened in 1919. Johann Radon was selected, but unfortunately, he accepted an offer at Hamburg University. The second candidate for the professorship in Brno, Roland Weitzenb¨ock, had already been appointed at the Prague German Technical University. Negotiations continued until 1921, when the rector of the Brno Technical University suggested the name of Horst von Sanden, who was a professor of mathematics working at the Clausthal Mining Academy. But Sanden rejected the offer, as did Robert König (1885–1979) and Georg Prange (1885–1941) in 1922. Eventually Karl Mayr, an assistant of mathematics in Brno before WWI, habilitated at the Vienna Technical University and became Tietze's successor in Brno. But this lasted only a short time. Dissatisfied with his position as extraordinary professor in Brno, Mayr left for the Technical University in Graz in 1924, even though he was offered no promotion [Šišma 2002, pp. 216–219]. Then, the year after Mayr left, the other mathematics professor at the Brno German Technical University, Schrutka, also quit to accept a professorship in Vienna. Negotiations to find proper replacements for Mayr and Schrutka were again long and difficult. In 1925, Rudolf Weyrich (1894–1971), who had been a student at Breslau University and a privatdocent in Marburg, was appointed extraordinary professor and in 1927, Lothar Koschmieder (1890–1974) came to Brno as an ordinary professor[Šišma 2002, pp. 220–225].

A further delaying factor in the negotiations was the fact that the Czechoslovakian Government often rejected foreign candidates. But graduates from the German universities in Prague or Brno were not sufficiently numerous to compete with Austrians and Germans. Only in exceptional cases (e.g., when a native expert could not be found) would the Ministry of Education consider a foreign candidate. Salary was another issue, as the Czechoslovakian crown's exchange rate was very low. As Josef Krames (1897–1986), who was professor of descriptive geometry at the German Technical University from 1927 to 1929, wrote to the Czechoslovakian Ministry of Education in 1927, the salary offered to him in Brno was that of an assistant in Vienna [Czech Nat. Arch., Josef Krames's personal file]. He nevertheless accepted the position in the hope that the status of extraordinary professor in Brno would help him get the professorship of descriptive geometry at the Graz Technical University. This indeed happened in 1929.

---

[24] "Co to znamená? To znamená, že akademický senát německé university v Praze neuznal právní stav tohoto státu, poněvadž připouští možnost, že by ta universita mohla býti přenesena do státu cizího, že tedy německá universita v Praze není, abych tak řekl, majetkem nebo statkem Československé republiky"(quoted from Czechoslovakian Parliament Debates, 34th Session, February 27, 1919).

Under these circumstances, the German Technical University in Brno was hardly in position to maintain its previous level of mathematical activities. Only fundamental mathematical lectures for engineers were now delivered; there were no more special lectures by privatdocents of mathematics. Mathematics professors kept their contacts with Austrian and German mathematicians and did all they could to get positions at foreign universities. They regularly participated in the meetings of the DMV and even organized such a meeting in Prague in 1929. But if we know of individual contacts between Czech and German mathematicians in Prague, such contacts seem to have remained extremely rare in Brno.

3.3. The Foundation of Masaryk University. Rather than confronting the German community by reducing the number of their institutions of higher education, the Czechoslovakian Government rather chose to develop new ones for their Czech nationals. For the Government, the support of existing institutions and the creation of new ones was a question of vital necessity, especially in the south of the country where the traditional road to Vienna for Moravian and Slovak students was now barred. In Brno, this political resolve translated into the enlargement of the Czech Technical University and the establishment of the brand-new Masaryk University.

In the Czech Technical University, a Faculty of Architecture was opened in 1919, after many years of requests. When peace returned, the number of regular students at the Czech Technical University jumped to more than 900. Technical universities were moreover opened to women after the war.[25] Due to the presence of the older classes, the age distribution changed slightly: in 1914 only 14% of the students were over 24 years of age; in 1919, they were 28%. Many of them were therefore in need of rapid qualification to start working as technicians for the new state. For this task, mathematics did not seem a priority. Indeed, despite the enlargement of the Technical University, the number of professors of mathematics, descriptive geometry, and mechanics remained the same after WWI. After the number theorist Lerch whose lectures at the Technical University had often been criticized for their high abstraction was appointed at Masaryk University (see below), Karel Čupr (1883–1956) was appointed as his successor in 1923. Descriptive geometry was taught by Pelíšek until 1928. In 1924 the only Slovak mathematics professor at a Czechoslovakian university, Jur Hronec (1881–1959), joined him.[26]

---

[25] Before 1919, women were only allowed to study at technical universities as exceptional students. Although there were several attempts to make technical universities open to women before the war, the Austrian Ministry of Education was always opposed to the idea. Women had been admitted as full students from 1897 on in the Empire universities but only in the Philosophy section. A tiny number of women were subsequently admitted to follow different studies. An interesting example of a woman scientist is the chemist Lise Meitner (1878–1968) who was admitted to university in 1901 and received a Ph.D. in Physics in 1905. On this question, the interested reader may refer to [Heindl & Tichy 1990].

Even before the collapse of the monarchy, Czech and Slovak leaders had already agreed on the necessity of creating two new universities, in Brno and in Bratislava. No later than during the second session of the Czechoslovakian Parliament on November 15, 1918, a group of deputies proposed the creation of a university in Brno with three faculties: philosophy, law, and medicine. Masaryk University was founded on January 28, 1919 with four faculties, the philosophical faculty being divided in two—philosophy and science—as had been done at the Prague universities. Before the beginning of the first academic year, professors were appointed in all faculties. Lectures started in the Faculty of Law and in the Faculty of Medicine in the academic year 1919–1920. The faculties of Science and Philosophy started their activities the following year. The main difficulty faced by the new university was the absence of rooms and buildings. All the faculties began their teaching in temporary quarters. Rooms were rented from the Technical University or other organizations. The first mathematical lectures by Lerch were held in a room where he had already taught his students of the Technical University. A huge plan for developing an academic quarter near the Technical University was not realized, and in fact only one building was built—the Faculty of Law.[27]

In the new Faculty of Science, there were two professorships for mathematics and two for physics. Lerch taught mathematical analysis and algebra, while Ladislav Seifert (1883–1956) was appointed professor of geometry.[28] Fréchet's correspondent Hostinský was appointed professor of theoretical physics while he was a privatdocent at the Prague University. Finally, Bedřich Macků (1879–1929), who had been extraordinary professor of physics at the Czech Technical University in Brno, was appointed professor of experimental physics. Those were distinguished recruits. As we have seen, Hostinský started his career with a mathematical thesis on Lie spherical geometry defended in 1906. The year he spent in Paris (1908–1909) was decisive for his scientific evolution and allowed him to prepare his habilitation which he defended on November 16, 1911 under the title "On Geometric Methods in the Theory of Functions." Back in Prague, Hostinský worked in secondary education, before he was called as privatdocent to the Prague University in 1912. In parallel with his secondary teaching, he began to give lectures on several topics of advanced mathematics (analytic function theory, differential geometry of curves and surfaces, differential equations, geometric applications of differential equations, etc.). Hostinský did not serve in the army during the war (most probably for health reasons) and stayed in Prague, where he was first acquainted with probability theory [Havlová et al. 2005]. Some months before his appointment in Brno, during the academic year 1919–1920, Hostinský taught Volterra's integral equation theory and their applications [Beránek 1984].

---

[26] Hronec worked in Brno until 1938. The fact that he was the only Slovak in this position illustrates the disparity between the two founding components of the multinational state. In 1939 Hronec was appointed professor and rector at the first Slovak Technical University in Košice [Franěk 1969, pp. 233–240].

[27] The other faculties were located in different parts of the town. The situation has not significantly changed until a new campus was inaugurated in September 2010.

[28] Ladislav Seifert is not to be confused with the German topologist Herbert Seifert (1907– 1996).

Seifert had also been privatdocent in Prague before he was appointed professor of geometry in Brno. In the academic year 1907–1908 he had studied in Strasbourg and Göttingen. An author of works devoted to algebraic geometry and differential geometry, he specialized in surfaces of the third order and quadrics of revolution. In differential geometry, he studied the properties of curves and surfaces and interpreted his results in descriptive geometry. An heir to the Czech geometers of the second half of 19th century, he remained outside the main trends in 20th-century geometry.[29] Seifert passed a double habilitation in 1920, at the Prague University in the field of mathematics and at the Technical University in descriptive and synthetic geometry. One may think that Masaryk University's urgent need for a professor of geometry accounted for this choice.

As we have seen, the creation of Masaryk University allowed Lerch to transfer and finally to lecture on mathematics, mostly to secondary school teachers. In 1920, Lerch shared with Hostinský the teaching of all mathematical and physical lectures and officially launched the mathematics department. Unfortunately, as he had become ill, Lerch restricted his activities mostly to teaching. His main contribution may have been to notice the talent of his student, Otakar Borůvka (1899–1995), whom he managed to appoint as the first assistant of mathematics at Masaryk University. Borůvka later became the leader of Brno mathematics in the second half of the 20th century.[30]

After Lerch's death in August 1922, Seifert became the head of the department and Eduard Čech (1893–1960) was appointed in 1923, only one year after his habilitation. Cech had studied mathematics and descriptive geometry at the Prague University from 1912. He was drafted during the war, but since he worked as a clerk in the rear for the next three years, he was able to read mathematical books and learn foreign languages (Italian, German, and Russian). After the war, he started to teach in a secondary school in Prague. He defended a thesis on differential geometry in 1920 and spent the academic year 1921–1922 in Turin with Guido Fubini. A geometer, he nonetheless was asked to teach mathematical analysis and algebra at Brno. After he was appointed full professor in 1928, Čech developed, as is well known, a strong interest for topology, a field in which he became one of the world's experts in the 1930s.[31]

Among the efforts made by Czech mathematicians at Masaryk University to develop research in Brno, let us underscore Hostinský's zealous direction of the library commission of the Faculty of Science, where he managed to obtain the inheritance of several personal book collections. Hostinský also devoted much energy to the founding of a journal Spisy vydávané přírodovědeckou fakultou Masarykovy university (Papers Presented to the Faculty of Science of the Masaryk University). He managed to conclude exchange agreements with hundreds of scientific institutions worldwide—a decisive point for the new university as, due to the low value of the Czech crown, foreign journals were unaffordable. Thanks to thesedeals, in 1937 the department of mathematics bought only six journals, and received around one hundred more through exchanges [Čech 1937]. Moreover, the creation of this new journal became a good opportunity not only for the faculty members but also for all Czechoslovak scientists to disseminate their results. In 1925, about fifty issues had been published as separate numbers.

---

[29] On Seifert, see [Hrdličková 2001]. About the Czech geometric school, see [Folta 1996].

[30] Borůvka was to become involved in a tight Czech-French collaboration. Borůvka's famous Spanning Tree algorithms were presented in Paris in 1926–1927 where he went to work with Elie Cartan on the suggestion of Hostinský and Čech. MU Archives contain several letters sent by Borůvka to Hostinský from Paris where he described how he had been welcomed by the Cartans. See [Nešetřil et al. 2001].

[31] On Čech, see [Gray 1994], [Katětov & Simon 1993].

## 4- A Franco-Czech Mathematical Axis : Dreams and Reality

Once Masaryk University was running, it appeared necessary for it to participate in international scientific programs and collaborations. More than an academic need, asserting a Czechoslovakian presence on the international scientific scene was an important political symbol. In mathematical science, a very interesting example of efforts to establish scientific relationships congenial to the new political map of Europe were, as mentioned above, the ties knitted between Fréchet in Strasbourg and Hostinský in Brno. The case is challenging because it combines several very different aspects: theaforementioned opportunistical consideration of the political situation of Europe in 1919but also the digestion and circulation of some ideas about chance and probability, in particular of Henri Poincaré (1854 -1912).

4.1 : Strasbourg as a laboratory for future universities

In December 1918, the French arrived in Strasbourg on a land they considered as a conquerred territory. They quickly became disillusioned and officials sent alarmist reports to Paris saying that in the newly returnedAlsace and Moselle, it was impossible to apply overnight laws and regulations governing institutions (including universities) in "inner" France. The first reason was linguistic : the students in their vast majority simply did not speak French. Moreover, all those who had already reached the age for higher studies, had started them in the German system.

Nevertheless, the French Government wanted the university to become a showcase of French science's successes [Olivier-Utard 2005, pp. 154–155]. The following letter from a French deputy to the administrator of Alsace, dated March 1919, is a good illustration.

> You know better than anyone the considerable importance that the Germans had given to this university and the coquettishness they showed to make of it one of the most brilliant, if not the most brilliant, of the whole empire. . . . When they left, they predicted that in less than three years France would have destroyed their work. How can we respond to this challenge?[32]

---

[32] "Vous savez mieux que personne l'importance considérable que les allemands avaient donnée à cette université et la coquetterie qu'ils ont mise à en faire une des plus brillantes sinon la plus brillante de l'empire. Vous avez certainement vu aussi qu'ils ont prédit en partant qu'en moins de 3 ans la France aurait saboté leur œuvre. Comment relever ce défi ?" Member of Parliament Maurice Maunoury to Alexandre Millerand, April 5, 1919 [Arch. Bas-Rhin, 1045W165].

Every effort was therefore geared toward turning Strasbourg into a first-rank French university. One came quickly to the conclusion that an original organization must take place in Alsace-Moselle. The aim was to "Frenchify" the university, and more generally the education system, while trying to keep as much as possible of the old German system. During the 1920s, Strasbourg became therefore an impressive laboratory of teaching and research initiatives where academics sought to invent the university of the future, with a freedom of action that seems to have had few equivalents in France. It is also noted that the officials from Paris, including the *Commissaire Général* Alexandre Millerand, became aware of the limited influence they had on the local situation. This lead them several times to leave local authorities handle with the problems. An exampleis the creation of the Institut d'Enseignement Commercial Supérieur (IECS) by the Strasbourg Chamber of Commerce. Reading the minutes of the 1919 meetings of that chamber shows how it had free range to create the new institute according its own criteria and agenda, with the central argument that the institution was necessary to provide Alsace with specialized technicians of trade and management adapted to the particularities of the local economy. Great freedom seems especially to have been left to the Chamber of Commerce to draw teachers for the new institution from the pool of the university professors.

4.2 International relations

To analyze thequestion of the international relations between Strasbourg university and Czechoslovakia, one must take into account French eagerness to create and tighten its relationship with Central and Eastearn European countries. Czechoslovakian leaders, such as Masaryk and Beneš, had strong personal and intellectual connections with France. French authorities devoted much effort to convince the Czechoslovakian Government and local administrations (universities, schools, cultural associations) of the importance of cultural and educational relations. Both universities recently established in Brno and Bratislava were the objects of special attention. After the war, it seemed clear that victorious France—and not Germany, its traditional competitor in this domain—would attract most advanced students from Czechoslovakia. While French officials betrayed their amazing self-confidence on that subject, German academics expressed their fear of French cultural hegemony. As early as in 1918, the historian and librarian Ferdinand Rieser (1874–1944) wrote that, if the Germans do not provide the necessary efforts to attract foreign students, "afterwards, the Russian and the Japanese will go to French universities which are not worse than German ones and will go back home and spread the French spirit."[33]

As said before, the French Government hoped to attract many foreign students from Eastern and Central Europe to Strasbourg.Soon, a rhetorical comparison between Czechoslovakia and Alsace appeared. Whenthe slavist André Mazon (1881-1967), who became professor at Strasbourg in 1919, invited the editor-in-chief of the journal *La parole tchécoslovaque à l'étranger* (Czechoslovak word abroad), Emil Čenkov, for the inauguration of Strasbourg university, the latter gladly accepted the invitation with the words : "there is so much similarity between Alsace-Lorraine and Bohemia !"[34].Both regions had just been rescued from what was perceived as the jaws of German imperialism. This idea influenced French cooperation with the new state in the political, economic, cultural, and, in particular, academic domains. Almost immediately after the war, the takeover of Strasbourg University and its reconstruction along French standards appeared as urgent tasks for the Government.[35]

---

[33] Akademische Rundschau, vol.5, 1918, p.322.

However, contrary to French authorities' first hunch, Strasbourg's strong ties with German culture was not an obvious advantage in the eyes of foreign students from Central Europe. Sent to Transylvania (a part of the Austro-Hungarian Empire assigned to Rumania) to attract students to Strasbourg in June 1920, a French diplomat reported the feeling of a local academic officials:

> For many people here, even among the most francophile, the Alsatian remains a hybrid being, as German as he is French, and who, condemned to be on one of the two sides, prefers France in the end where life is sweeter. I felt that [they feared] that the air one breathes in Strasbourg is not "absolutely pure French air." You will certainly recognize therein the effect of the Fritz propaganda on the mind of these brave Transylvanians.[36]

4.3 Fréchet and mathematics in Strasbourg

Before being appointed to Strasbourg as professor of higher analysisin 1919, Fréchet had been professor at the University of Poitiers between 1910 and 1914. He was already famous in the international mathematical community due to his outstanding 1906 doctoral thesis on the topology of functional spaces[Taylor 1982]. In addition to his mathematical fame, Fréchet had another asset: he was a polyglot. A singular fact of his life is the energy he devoted to the promotion of Esperanto, a language in which he wrote several mathematical papers.[37] His knowledge of English was excellent at a time when this was unusual. During the war, he had served as an interpreter for the British Army. And he knew German well, a useful asset in Strasbourg.

Fréchet's speech at the official opening of the university in November 1919 was remarkable. Far from subjecting uniformly all that the Germans had made in

---

[34] Il y a tant de ressemblance entre l'Alsace-Lorraine et la Bohême. (letter from Emil Čenkov to André Mazon, 18 October 1919).

[35] On the history of Strasbourg University as stakes of power between France and Germany, the interested reader may refer to [Crawford & Olff-Nathan 2005].

[36] "Pour beaucoup de gens ici, et des plus francophiles, l'Alsacien demeure un être hybride, autant allemand que français, et qui, réduit à vivre dans l'un des deux camps, préfère en définitive celui de la France où on lui fait la vie plus douce. J'ai senti [qu'ils craignent] qu'on ne respire pas à Strasbourg 'un air de France absolument pur'. Vous reconnaissez là, bien certainement, Monsieur le Doyen, les effets de la propagande boche sur l'esprit de ces braves transylvains." Letter from the military attaché in Bucarest to the Dean of Strasbourg University, June 6, 1920 [Foreign Affairs, file France, 43, F◦ 128].

[37] See for example Fréchet's report on the international conference for the use of Esperanto in science, [Fréchet 1926].

Strasbourg to public scorn, he made much of the fact that the German higher education system had some advantages it would be good to keep, such as the contact of students with scientific research through Seminars, which occurs much earlier than in France. This balanced speech (even if it is stained here and there with unavoidable patriotic flights), so soon after the Armistice, clearly departed from the ambient ostracism of the time against anything with a German flavor. This is probably why Borel, during the brief attempt to revivehis journal*Revue du Mois* after the war, decided to publish Fréchet's speech in 1920 as [Fréchet 1920].

Several French scientists active in the organization of wartime scientific mobilization, especially at the *Direction des inventions intéressant la défense nationale* (Board of Inventions) where Emile Borel played a major role, were chosen to form the new faculty of Strasbourg University. Among them, some understood the importance of developing statistics. During the imperial period, Strasbourg had already, with Wilhelm Lexis and Georg Friedrich Knapp, been a major center for the discipline.[38]

When the IECS was created, Fréchet and the sociologist Maurice Halbwachs (1877-1945) were among the first teachers. Halbwachs is a very singular personality of French intellectual scene of the time. While lying in the tradition of the sociology of Durkheim, he had learned with Lexis during a stay at Göttingen how the use of mathematical statistics played a major role in modern sociology. His *thèse complémentaire*(complementary thesis) published in 1912 ([Halbwachs 2010])examined the work of Quetelet and moral statistics. It reveals how the author was a fine reader of recent works on probability and on the interpretation and analysis of the presence of randomness at the heart of science, particularly those of Poincaré[39].

The meeting between the two Maurice at the Institute of Strasbourg and their respective courses on statistics and Insurance resulted in an original book, published in 1924 as [Fréchet & Halbwachs 1924][40].In the preface, the authors explain how the dialogue between a mathematician and a sociologist proved fruitful. The book participated in the renewal of mathematical tools in the humanities and is closely connected to the dissemination of the concepts of mathematical statistics in the 1920s in France studied in [Catellier and Mazliak, 2012].

At the same time he was involved in these Strasbourg experiments, immediately after his nomination as head of the mathematical department, Fréchet began a campaign to promote Strasbourg mathematics abroad.Under these circumstances, it seems understandable that Fréchet would take the initiative of writing his Czech colleagues to receive recent information about the state of mathematics in their country.[41]

---

[38] On the development of mathematical statistics in France, see [Catellier & Mazliak 2012]. See also [Armatte 2001].
[39] We can not dwell here on Halbwachs but the interested reader is invited to refer to [Baudelot Jaisson 2007] or [Brian Jaisson 2005].
[40] On the genesis and the contents of the book see [Brian&Mazliak 2015].
[41] On Fréchet's effort to organize mathematics in Strasbourg and, especially, the international dimension of this effort, see also [Siegmund-Schultze 2005, pp. 187–189].

As we have seen, the letter he wrote on June 29, 1919, was finally answered by Hostinský in October. On November 12, 1919, Fréchet answered Hostinský's letter. He suggested that French abstracts to Czech publications be collated in a single journal. With a tinge of paternalism, Fréchet cautiously wrote that he had thought it "interesting to let [Hostinský] know the opinion of a stranger who wishes nothing but good things for Czech scientists and mathematical science," and that collecting these abstracts would show that there was a large part "due to Czech science in what has usually been attributed to the Germans in Austria."

In a letter dated June 1, 1920, written on official letterhead of the organization committee of the sixth International Congress of Mathematicians (ICM, to be held in Strasbourg in September of that year), Fréchet enclosed a little brochure entitled "The Teaching of Mathematics at the University of Strasbourg" and printed to attract students to the Alsatian capital. He asked Hostinský to publish it in a Czech journal. According to Fréchet, the Strasbourg University needed to stimulate the influx of new students, so propaganda was needed for a little while.[42] In response at the end of June, Hostinský wrote that he hoped to meet Fréchet for the first time at the ICM in Strasbourg. The Czech delegation to the Congress was important (11 members out of a total of some 200 persons). After the congress, the mathematician Bohumil Bydžovský (1880–1969), who was one of the most active members of the Czech mathematical society and had been appointed the same year as a tenured professor at Prague University, reported that

> the contact with mathematicians from the University of Strasbourg, our main partner in the West, was particularly cordial. The interest they showed for our scientific, educational and social situation seems to warrant that reciprocal exchanges be pursued, obviously for the prosperity of our science [Bydžovský 1921, p. 47].

4.4 Hostinský and France

We have already mentioned the study period of Hostinsky in Paris, where he encountered French leading mathematicians and studied their work. Another element played a direct role in the story presented here, namely his discovery of Poincaré's work on probability and randomness. As we mentioned earlier in the case of Maurice Halbwachs, the circulation of Poincaré's texts at the beginning of 20th century is an important point of the scientific life of the time, as Poincaré played a central role to bring scientists in contact with a new physics involving a use of probability theory.[43] However, as well expressed Bernard Bru ([Bru, 2003]) one easily gets the impression that everyone at the time read Poincaré's philosophical texts (especially in *Science and Hypothesis* or *Science and Method*) but that few really understood him. And in fact, as During mentions ([During 2001, p.93]), Poincaré leaves the question of how his considerations implied a redefinition of the notion of truth to the care of logicians and metaphysicians (and, as we have seen in Halbwachs' case, of sociologists). Exceptions were rare who studied Poincaré's probabilistic considerations in depth. We have already mentioned Borel and Halbwachs. Now we shall turn to an outsider, coming from nowhere, the Czech philosopher Karel Vorovka.

---

[42] As a consequence of this letter, we only were able to find a brief publication of the list of lectures offered in Strasbourg university in the journal Casopis (Journal for the Cultivation of Mathematics and Physics) in the Fall 1920. Compared to the rather luxurious presentation of the brochure sent by Fréchet, the Czech advertisement seems a bit dull.

[43] See [Mazliak 2014] for a detailed study of Poincaré and probabilities.

An interesting and very complete study on Vorovka was published in Czech some years ago ([Pavlincova, 2010]) and hopefully it will become more accessible in a more widely known language. Some complements can also be found in [Mazliak 2007b] and the references therein. Two reasons explain this general ignorance of Vorovka: the fact that his works, mostly in Czech, were never translated, and also that his early death precluded the collection of his ideas in a large-scale work. Placing himself in the tradition of Bernhard Bolzano (1771-1848), the major figure of the philosophical scene in Prague during the 19th century, Vorovka sought an approach of the question of truth by combining both his strong mathematical education and a rather strict religious philosophy, an original syncretism of empiricism and idealism which had close links with the works of the hero of the Czechoslovak independence, the philosopher T.G. Masaryk, and with American pragmatist philosophy, in which he had much been interested (see [Pavlincova, 2010, pp.67-69]).

Vorovka's discovery of Poincaré's philosophical writings at the beginning of the 20th century was a real revelation: he drew from them the conviction that the scientific discoveries at the end of the 19th century, especially in physics, compelled a reconsideration of the question of truth and free will. Vorovka showed a real originality in that he did not content himself with principles, but closely studied the mathematical problems raised by the theory of probability. He was a diligent reader of Bertrand's textbook, of Borel's texts, but also of Markov's works, publishing several works inspired by papers of the Russian mathematician. At the time when he was granted tenure at the Czech University in Prague, around the year 1910, Vorovka metBohuslav Hostinský, who had just returned to Bohemia after his period of research in Paris. In Hostinský's own words ([Hostinský 1929]), it is through the discussions he had with Vorovka that he learned about Poincaré's works, and he began to reflect upon the calculus of probability. In 1917, Hostinskýpublished his first paper on probability in the *Proceedings* of the Czech Academy [Hostinský 1917], devoted to a new solution of the Buffon needle problem. This famous problem of the calculus of probability involves the determination of the probability that a needle thrown over a parquet floor intersects one of the regularly spaced grooves. Its classical treatment of this problem is based on the hypothesis of uniform distribution of the (random) angle between the needle and the grooves' direction and of the (random) distance between the middle of the needle and the closest groove. It amounts to consider that the floor is infinite in every direction and that the needle could be thrown anywhere on this surface. Hostinský criticized these assumptions, which he found absurd physically, and looked for the result within the reasonable framework of a finitely extended floor. To achieve this, he made use of Poincaré's arbitrary function method and showed that for any arbitrarily given distribution of the random angle and distance, if the number of grooves goes to infinity, the value of the desired probability converges to the one obtained under the illegitimate classical assuptions.

4.5 Strasbourg conference and its consequences

At the Strasbourg ICM where he indeed met Fréchet, Hostinský gave two talks, one on differential geometry and the other on mechanics. During the spring of 1920, Hostinský had sent Picard the translation of his 1917 paper. As soon as he received it on April 18, 1920, Picard offered to publish it in the *Bulletin des sciences mathématiques* [Hostinský 1920]. In a letter to Hostinský dated November 7, 1920, Fréchet wrote that he had carefully read this slightly modified version of the 1917 paper. From this reading, Fréchet was motivated to write his very first article in the field of probability theory [Fréchet 1921]. At the same time, the interaction was, for Hostinský, the occasion for making a first step in the direction of the ergodic principle to which he would later devote important studies. Fréchet had been drawn to probability theory for wider reasons than his nascent correspondence with Hostinský. For one thing, Fréchet's joint lectures at the IECS, and more generally his teaching duties at university, certainly played an important role. However, it is clear that his continuing correpondence with his Czech colleague was a major incentive for his interest in Markov chains since, in his letters, Hostinský discussed all new developments on this topic. Fréchet himself wrote after Hostinský's death that:

> Among his various topics of research, [Hostinský] drew my attention to the theory of probabilities in chain. Hence, it is thanks to him that I could write the second volume of my studies on modern probability theory, and in the book I frequently used his ingenious methods.[44]

Hostinský's contacts with French mathematicians went beyond his correspondence with Fréchet. In a note submitted to the Paris Academy of Sciences in 1928, Hostinský introduced an elementary version of the ergodic theorem for a Markov chain with continuous state [Hostinský 1928]. Hostinský's work on this topic came before the spectacular development of the 1930s at the hands of Andrei N. Kolmogorov and others. Upon reading Hostinský's article, Jacques Hadamard plunged into probability for the first and only time of his life, a period referred to as his "ergodic spring" which ended at the Bologna ICM in September 1928 where Hadamard gave a talk on the ergodic principle [Bru 2003, pp. 158–159]. Between February and June 1928, Hostinský and Hadamard exchanged many letters, published several notes responding to one another, and also met during Hadamard's journey to Czechoslovakia in May. From this moment, Hostinský acquired real international prestige, and in the 1930s, his little school in Brno became an active research center on Markovian phenomena.

Hostinský's glorious year in 1928 however marked the end of the special relations between Brno and Strasbourg. That year, Fréchet went back to Paris to help Borel with the new Institut Henri Poincaré. This was the institution where Hostinský's international contacts were to develop in the 1930s. In their numerous exchanges during these years, Fréchet and Hostinský tackled many aspects of Markov chains with discrete and continuous states. This tends to prove that Hostinský and Fréchet had become each other's main reader on the topic during the 1930s. In 1936, Fréchet advised Wolfgang [Vincent] Döblin to write to Hostinský, and this resulted in a short but interesting correspondence studied elsewhere [Mazliak 2007]. Hostinský was moreover eager to lean on the prestige acquired from his French relations for further development of international contacts for mathematics in Brno. During the 1930s, his activity was particularly intense to federate the mathematical research of the Central and East European countries. He was one of the main organizers of the first two congresses of mathematicians from Slavic countries (in Warsaw in 1929, Prague in 1934). His personal archive in Brno testifies to his energy, especially for inviting the Soviet mathematicians, which mostly failed due to the conditions of the Stalinist era in the USSR. Until its annihilation in the dramatic events of the German annexation of 1939 and World War II, the Brno school of mathematical physics remained one of the most successful mathematical centers in the new Central Europe that came out of WWI.

---

[44] "Parmi toutes ses recherches si variées, il avait su m'intéresser à la théorie des probabilités en chaîne. De sorte que c'est grâce a lui que j'ai été amené à écrire sur ce sujet le second livre de mes recherches sur la théorie moderne des probabilités, ouvrage où j'ai eu à invoquer ses ingénieuses méthodes en de nombreux passages" [MU Archives, Hostinský's personal file; Fréchet to the Rector of Masaryk University, May 5, 1951]. Fréchet alludes to the book published in 1937 as part of Borel's Traité du calcul des probabilités et de ses applications (Gauthier-Villars, Paris) under the title Méthodes des fonctions arbitraires. Théorie des événements en chaîne dans le cas d'un nombre fini d'états possibles.

**Conclusion**

The creation of Czechoslovakia in 1918 was an attempt at reorganizing European states after the end of World War I. In places where there was a tense cohabitation of several national communities (as in many parts of the collapsed Austro-Hungarian empire), it was necessary to choose a form of organization that allowed the coexistence of several traditions. This was in particular the case with the organization of the educational system.

When one studies the local case of Brno, the capital of Moravia on the border of Austria, it is vital to understand how the difficult contacts between the Czech majority and the large German minority influenced the edification of educational institutions between 1880 and 1930. Though the German minority lost its domination in 1918, the institutions were still much influenced by the culture that had prevailed before the war. Despite the attempts made to stir interest in favor of the Allied countries, and despite the significant political changes that occurred in the region, unsolved contradictions, antagonisms, and lack of communication characterized the relations between France and Czechoslovakia. French eagerness to become the strongest supporter of the young Czechoslovakian state gave rise to misunderstanding and suspicion on the part of Czechoslovakians. After promising beginnings illustrated by a proclaimed proximity between Strasbourg, Prague, and Brno, the relation was in fact never developed to the extent the French had hoped.

As for mathematics, the discipline was mainly active in the German Technical University before the war with the creation of a local German Mathematical Society when Hamel was given a professorship in mechanics, and in the new Masaryk University after the war where Hostinský and his international contacts, especially with France, played and important role. Hostinský thus succeeded in federating a very active kernel of research in mathematical physics in Brno, mainly studying the structure of Markovian processes. This was an important occasion for Brno to become the center of an international network of mathematicians, including in particular scientists from Eastern Europe. Significantly, this network was efficiently structured around the exchange of scientific publications by Hostinský. But this did not survive to the German occupation in 1938. The Second World War, followed by the new geopolitical regime of alliances with the Soviet Union during the Cold War, turned Brno mathematicians' gaze eastward.